\documentclass[10pt,letterpaper]{article}
\usepackage{llncsdoc}
\makeatletter
\DeclareRobustCommand{\text}{%
  \ifmmode\expandafter\text@\else\expandafter\mbox\fi}
\let\nfss@text\text
\def\text@#1{{\mathchoice
  {\textdef@\displaystyle\f@size{#1}}%
  {\textdef@\textstyle\f@size{#1}}%
  {\textdef@\textstyle\sf@size{#1}}%
  {\textdef@\textstyle \ssf@size{#1}}%
  \check@mathfonts
  }
}
\def\textdef@#1#2#3{\hbox{{%
                    \everymath{#1}%
                    \let\f@size#2\selectfont
                    #3}}}
\makeatother
\usepackage{amssymb}
\usepackage{graphicx}
\usepackage{epsfig}
\usepackage{wrapfig}
\usepackage{epstopdf}
\usepackage{algorithm}
\usepackage{algpseudocode}
\usepackage{url}
\usepackage[toc,page]{appendix}
\usepackage[letterpaper,top=3.6cm,bottom=4.5cm,hmargin={3.45cm,3.45cm}]{geometry}
\begin{document}

\title{Physics driven real-time blood flow models for virtual interventional planning}


\author{Sethuraman Sankaran, David Lesage, Rhea Tombropoulos, Nan Xiao, Hyun Jin Kim, David Spain, Michiel Schaap and Charles A. Taylor\tnoteref{cauth1}  } 


{\Large
\textbf\newline{Physics driven reduced order model for real time blood flow simulations}
}
\newline
\\
\begin{flushleft}
Sethuraman Sankaran\textsuperscript{1},
David Lesage\textsuperscript{1},
Rhea Tombropoulos\textsuperscript{1},
Nan Xiao\textsuperscript{1},
Hyun Jin Kim\textsuperscript{1},
David Spain\textsuperscript{1},
Michiel Schaap\textsuperscript{1},
Charles A. Taylor\textsuperscript{1,*}
\\
\bigskip
\bf{1} HeartFlow Inc., 1400 Seaport Blvd, Building B, Redwood City, CA 94063, U.S.A.
\\
\bigskip

\end{flushleft}

\begin{abstract}
Predictive modeling of blood flow and pressure have numerous applications ranging from non-invasive assessment of functional significance of disease to planning invasive procedures. While several such predictive modeling techniques have been proposed, their use in the clinic has been limited due in part to the significant time required to perform virtual interventions and compute the resultant changes in hemodynamic conditions. We propose a fast hemodynamic assessment method based on first constructing an exploration space of geometries, tailored to each patient, and subsequently building a physics driven reduced order model in this space. We demonstrate that this method can predict fractional flow reserve derived from coronary computed tomography angiography in response to changes to a patient-specific lumen geometry in real time while achieving high accuracy when compared to computational fluid dynamics simulations. We validated this method on over $1300$ patients that received a coronary CT scan and demonstrated a correlation coefficient of 0.98 with an error of $0.005 \pm 0.015$ (95\% CI: (-0.020, 0.031)) as compared to three-dimensional blood flow calculations.  
\end{abstract}
{\let\thefootnote\relax\footnote{{* Corresponding author: Email: taylor@heartflow.com, Ph.No: (650) 740 8222}}}

\section{Introduction}

Patient-specific modeling of blood flow has emerged as a tool with increasing importance in the diagnosis and treatment of patients with coronary artery disease~\cite{TaylorMin13}~\cite{HFNXT}~\cite{Douglas}~\cite{Patel}~\cite{Kim}~\cite{PACIFIC}. A unique strength of simulation methodologies lie in the predictive modeling of hemodynamics in response to unplanned events (such as progression or regression of lesions), or the outcome of planned procedures (such as surgical intervention)~\cite{TaylorPredictive}. For predictive modeling to be used in the cardiac catheterization laboratory or other invasive procedures, it is imperative that the modeling tools can generate results in seconds.\par

In this work, we address one of the most challenging aspects of interventional planning for patients with coronary artery disease - Which coronary artery stenoses are having the greatest impact on blood flow?  Due to the vastly different anatomies and lesion morphologies in a normal patient population, it is infeasible to have a preset approach for the treatment strategy. Therefore, the planning of invasive procedures in a patient is left to the knowledge, intuition and experience of the clinician based on available data. This paper describes a patient-specific framework that can virtually model different scenarios and help the physician evaluate various strategies by identifying and ranking the impact each coronary artery stenosis has on the blood flow. We describe a fast and accurate tool derived from a patient specific coronary CT angiography scan that is applicable in a clinical setting. The clinical metric computed is the fractional flow reserve (FFR) which is measured in the cardiac catheterization lab with a pressure wire inserted in the patient's coronary arteries. FFR is ratio of the time-averaged pressure downstream of a coronary artery stenosis to a time-averaged reference aortic pressure under conditions of maximum hyperemia typically induced by the intravenous administration of adenosine. FFR computed from coronary computed tomography data ($\text{FFR}_{\text{CT}}$) is derived by simulating blood flow and pressure in a patient-specific anatomic model using computational fluid dynamics (CFD) methods~\cite{TaylorMin13}.\par

Fast assessment of blood flow has been an active area of research in the past several years. Recently, many successful methods have evolved that apply simulation methodologies for the estimation of clinical quantities of interest from medical imaging data. For instance, artificial intelligence methods have been applied to quantify hemodynamics from phase contrast MRI scans~\cite{Arterys}, or to quantify information from handheld ultrasound devices~\cite{Butterfly}. While these techniques help in the automation, miniaturization or reproducibility in the extraction of information, the approach, fundamentally, is to quantify information already present in the medical images. In contrast, this paper focuses on predictive modeling of the effect of treatment, for which no current patient-specific non-invasive alternatives exist.\par

The development of reduced order models for the simulation of blood flow has been the focus of several studies~\cite{Mette1,Mette2,Quarteroni1,Quarteroni2,Quarteroni3,McLeod,Ohhara}. Kassab and colleagues~\cite{Kassab} proposed an analytical model based on  the conservation of energy, and considered various energy losses associated with a lumen narrowing such as convection and diffusion and also accounted for the losses associated with sudden constriction and expansion in lumen area. They demonstrated the performance on a tube-like stenosis model, however calculation of some of the higher order geometric quantities (higher order gradients of radius or area) accurately on patient-specific geometries is not possible due to the imaging resolution. Schrauwen et al.~\cite{ROM1}~\cite{ROM2} improved upon the traditional reduced order models that assume parabolic velocity by deriving a velocity profile based on the geometry and flow, and validated it on straightened coronary arteries. However, the velocity profiles on patient specific models with multiple bifurcations and interfering stenosis (effect of one stenosis overlapping on another) can be significantly different. Nithiarasu and colleagues~\cite{Nithiarasu}~\cite{Nithiarasu2}~\cite{Nithiarasu3} developed reduced order models that accounted for intramyocardial pressure and material properties of the arterial wall, and evaluated performance in a virtual cohort of 30 lesions with the same global geometry. Itu et al.~\cite{Itu} proposed a machine learning approach for assessment of fractional flow reserve and verified it against another reduced order model on synthetic data. While these approaches have helped immensely in understanding how the complex physics of fluids can be captured with a reduced order system, these are not applicable for clinical translation due to the high demands on accuracy, latency and robustness. Our goal in this paper is to develop an algorithm that meets these requirements by demonstrating a very high accuracy over a large cohort of patient-specific models with diverse charecteristics. \par

The framework we propose uses a response surface~\cite{ResponseSurface} methodology and can predict the results of the simulations significantly faster and with an accuracy close to the high-fidelity simulations. The reduced order model is parameterized using the response surface built using a full order high fidelity high cost (HFHC) model performed at certain pre-set configurations. Simulations for the HFHC model may utilize all the information available about the system (such as using the full spatial and temporal representation). A response surface is a mathematical relationship between a quantity or quantities of interest or parameters, and an underlying representation or property, and have been built successfully for optimization problems~\cite{SankaranPoF}.  The pre-set configurations involve exploration of both the space of lumen geometries and flow rates. As we increase the number of HFHC simulations used to construct the response surface, the accuracy of the real-time predictive model becomes closer to the HFHC model within some error margin. The response surface is used to explore the parameteric space by using a reduced order model to interpolate the solutions. As a result of using physics driven interpolatory functions derived from one-dimensional Navier-Stokes equations, our approach can mimic the results of the computational  fluid dynamics simulations better than other standard approaches such as using polynomial or Lagrange interpolatory functions. The response surface is built to exactly match the output of the high-fidelity model for the configurations where high-fidelity simulations were performed.  \par

We demonstrate that we can achieve excellent accuracy of the proposed method compared to three-dimensional simulations for coronary blood flow using just four full order offline simulations. We also evaluated the performance of the algorithm in predicting hemodynamics post stenting or percutaneous coronary intervention (PCI) against invasive data for different lesion types. PCI is the process of inserting a catheter with a balloon attached to restore and attain ideal lumen geometry on a section of the blood vessel. We demonstrate that we can achieve high accuracy against invasive data and that it can be used as an interactive tool to plan PCI. \par

We structure the paper as follows. In section $2$, we describe the method behind the patient-specific reduced order model. We describe the HFHC model and demonstrate how we can build a real-time model using a set of HFHC models. In section $3$, we demonstrate the performance of the method against HFHC simulations. We use internal data that encompasses different lesion types and FFR ranges to validate the model. We finally discuss implications of the method and clinical translation in section $4$. \par

\section{Methods}
Blood flow in the human arterial system can be effectively modeled using computational fluid dynamics~\cite{TaylorMin13}~\cite{TaylorHughes}. These simulations can be performed on complex patient-specific anatomic geometries that are reconstructed from coronary CT scans. The Navier-Stokes equations are used to model blood flow and lumped parameter boundary conditions are prescribed that account for the flow-pressure relationship in the microvasculature (vessels that are not visible in the CCTA and therefore not modeled). Vignon-Clementel et al.~\cite{Clementel} developed methods for coupling boundary conditions using a circuit analogy of the microvasculature with the Navier-Stokes equations and Kim et al.~\cite{Kim1} specialized this method to modeling boundary conditions for coronary arteries. Sankaran et al. used a set of ordinary differential equations~\cite{SankaranAnnals} instead of an analytical formulation for modeling blood flow in the microvasculature. In this work, we use resistance (ratio of pressure drop to flow) boundary conditions for the microvascular network, which are now routinely used in clinical practice~\cite{TaylorMin13} for the non-invasive assessment of $\text{FFR}_\text{CT}$. \par

\subsection{Computational Fluid Dynamics Model} 
The partial differential equations that model blood flow are represented as  \\
\begin{equation}
                                                   \mathcal{N}(u;p)=0     \hspace{0.5in}       \text{in    } \Omega
\label{eqhfhc1}
\end{equation}
with the boundary conditions
\begin{equation}
                                                               b(u;p)=0           \hspace{0.5in}           \text{in    }  \Gamma
\label{eqhfhc2}
\end{equation}

where $\mathcal{N}$ is a non-linear differential operator modeling the Navier-Stokes equations, $u$ are the blood velocities and pressures, $p$ represents fixed parameters, $\Omega$ is the problem domain and $\Gamma$ is the boundary of the domain. Since this model is already validated~\cite{HFNXT}~\cite{PACIFIC} and used in clinical practice for the calculation of $\text{FFR}_{\text{CT}}$, we use the model governed by Equations~\ref{eqhfhc1} and~\ref{eqhfhc2} as our HFHC model. Our goal in this work is to build a low cost model that is equivalent in accuracy to the high fidelity model.\par

\subsection{Reduced Order Models}

A reduced order model of the partial differential equation approximates the operator, ${\mathcal N}$, using a simpler operator, $\hat{\mathcal N}$ (e.g., ordinary differential equations), reduces the dimensionality of the solution space $u$ and/or reduces the dimensionality, of the problem space $x$, to a lower dimensional space, $\hat{x}$, $\hat{u} \equiv \hat{u}(\hat{x} )$ and simplifies the parameter set $p$ to $\hat{p}$. In this case, $x$ can be the three dimensional representation of the geometry, and $\hat{x}$ can be a corresponding one-dimensional representation. The reduced order model problem is posed as follows,  where we omit the spatial dependence of $u$ for clarity: \\
\begin{equation}
                 {\hat {\mathcal N}} (\hat{u};  \hat{p})=0       \hspace{0.5in}              \text{in    }  \hat{\Omega}
\label{romeq1}
\end{equation}
with the boundary conditions \\
\begin{equation}
                                                               \hat{b}(\hat{u}; \hat{p})=0      \hspace{0.5in}            \text{in    }  \hat{\Gamma}
\label{romeq2}
\end{equation}
The goal is to have $\hat{u} (\hat{x})$ be a good approximation of $u(x)$ at $\hat{x}$. To achieve this, our approach involves performing simulations of the original HFHC model for various geometries, $\Omega$, boundary domains, $\Gamma$ and boundary conditions $b(u;p)$ such that we can generate a response surface and subsequently build a good approximating solution to the problem. In order to generate the response surface and to ensure that we limit the exploration space to realistic geometries, we can impose bounds on the domain of interest and the boundary conditions that the system will be subject to. The original governing equation is solved in a series of domains, boundary conditions and parameters. \\
\begin{equation}
\left( (b_1 (.),\Omega_1, \Gamma_1, p_1 \right), \left( (b_2 (.),\Omega_2, \Gamma_2, p_2 \right), \cdots,\left( (b_M (.),\Omega_M, \Gamma_M, p_M \right)
\label{Mconfig}
\end{equation}
where $M$ is the number of HFHC simulations performed. The results of the HFHC simulation for each of these configurations are
\begin{equation}
\left( u_1, u_2, \cdots u_M \right)
\end{equation}

\begin{figure}[h!]
\begin{center}
\includegraphics[width=12cm]{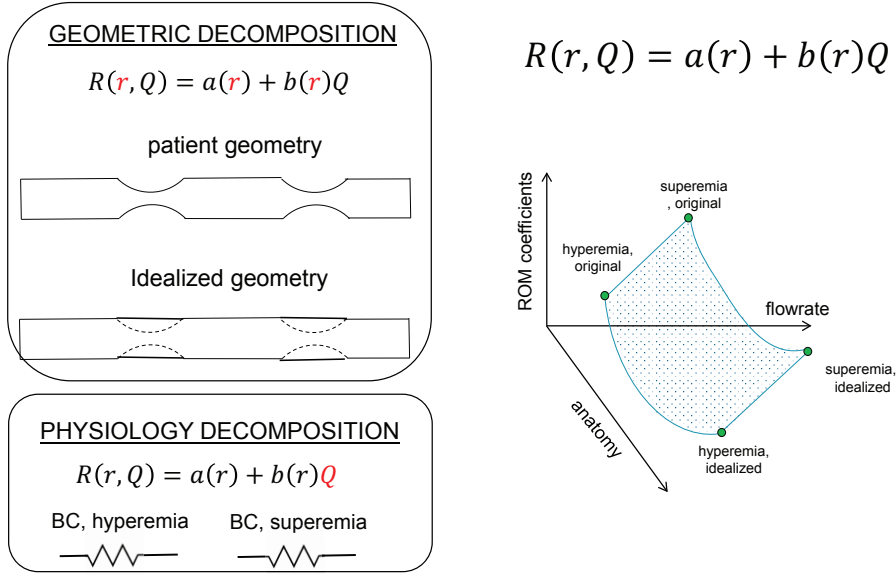}
\caption{(left) The reduced order model has two terms: a purely geometric term $a(r)$ (parameterized by radius $r$) and a term that varies linearly with flowrate, $Q$, with a geometry dependent slope, $b(r)$. The coefficients are fit for each patient based on the patient geometry and idealized geometry, and (right) four simulations are performed on two lumen segmentations (patient and idealized) with two boundary conditions provided for each.
\label{fig:twoterms}}
\end{center}
\end{figure}

\subsection{Physics driven response surface method}
Any sampling method or quadrature method can be used for the selection of the $M$ configurations in Equation~\ref{Mconfig}. As each simulation is computationally expensive, sampling methods such as Monte-Carlo or latin hypercube sampling~\cite{LHS_main} tend to converge slowly in building an accurate representation of the underlying function. Functional space methods, such as stochastic collocation or adaptive stochastic collocation methods~\cite{XiuCollocation}~\cite{SankaranSC}, can help accelerate the convergence for a class of functions. However, they still demand significant calculation to generate the response surface. We improve upon this further by regularizing the problem using a physics driven interpolation approach. Instead of using polynomials, which is commonly used for response surfaces, we use the analytical solution from one-dimensional, steady state Navier-Stokes equations to interpolate within the response surface. This enables employing a significantly reduced number of simulations while still enabling equivalence with the HFHC simulations. In fact, we later demonstrate that four HFHC simulations at the extremas of the exploration space are sufficient to obtain equivalence with the full order model (refer Fig.~\ref{fig:twoterms}). \par

A patient-specific response surface, ${\mathcal R}$, is a mapping of the coefficients of the resistance given a geometry, boundary conditions and parameters,  
\begin{equation}
\hat{c} \sim {\mathcal R}(\Omega, p, b(.), \Gamma)
\end{equation}
wherein $\hat{c}$ captures the complexity of the original equations, enabling $\hat{\mathcal N}$ to be a less complex operator than $\mathcal{N}$. The reduced order model, ${\mathcal R}$, can be constructed such that  $\hat{u}(\hat{x}) \equiv u(\hat{x})$ at the $M$ configurations where the HFHC simulations are performed. This approach lets us solve the reduced order problem (refer to equations~\ref{romeq1} and~\ref{romeq2}) for the unknowns $\hat{u}$, while ensuring that the results are identical at the $M$ configurations. The solution field $\hat{u}$ contains only flowrates and pressures instead of three velocity components and pressures. The reduced model is evaluated only along the vessel centerlines along the vessel path instead of the entire volumetric patient-specific geometry. The reduced order operator $\hat{\mathcal{N}}$ contains the affine resistance model that relates flowrates to pressures and the flow-split model at bifurcations, which are explained in later sections. Naturally, the approximations to the HFHC solutions at the intermediate configurations will be better for larger $M$, but so will the time needed for the offline computations. These are described in detail in the following sections. \par

\begin{figure}[h]
\begin{center}
\includegraphics[width=17cm]{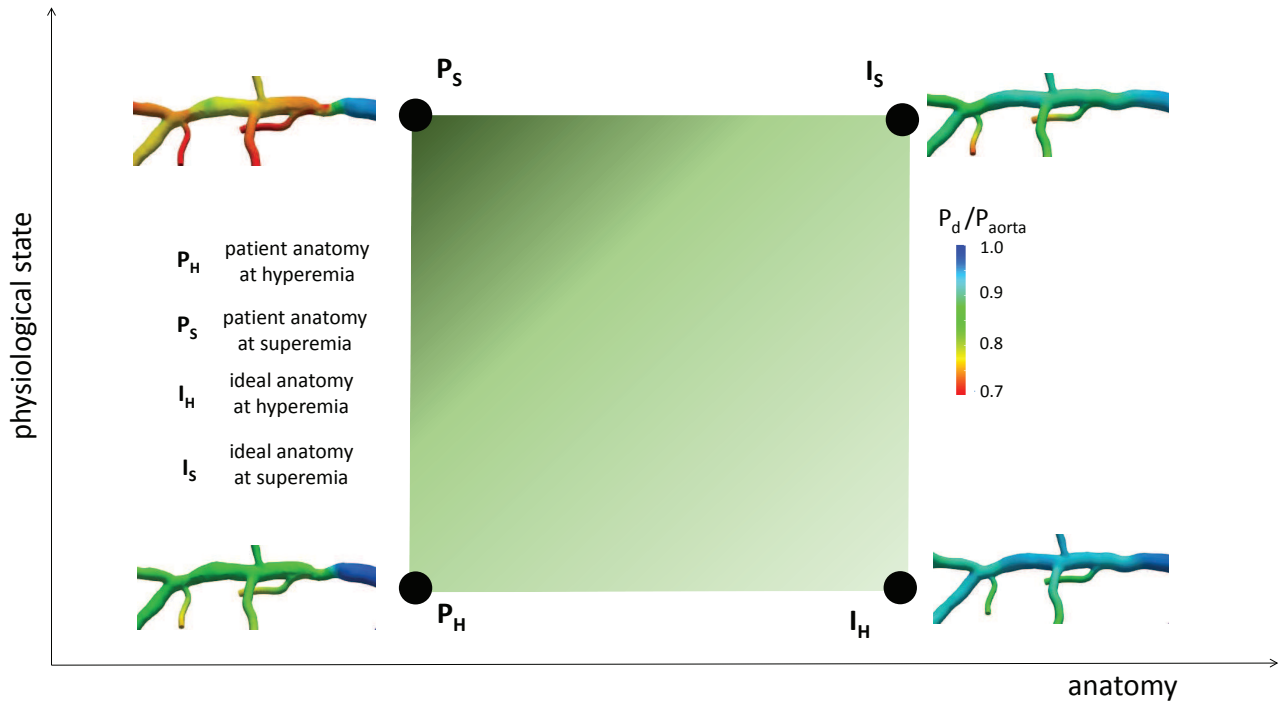}
\caption{Illustration of the proposed method that builds a patient-specific reduced order model which depends on both the anatomy and the physiologic state. Four simulations are used to explore the embedding space of the lumen geometry and physiologic state. Hyperemia and Superemia refer to different boundary conditions (or flowrates) and the colormap within the window represents how the reduced order model resistance varies within the patient-specific exploration space. $\text{P}_{\text{d}}$ is the local pressure and $\text{P}_{\text{aorta}}$ is the reference aortic pressure. \label{illustration} }
\end{center}
\end{figure}

\subsection{Patient-specific ideal geometries}

\begin{figure}[h!]
\begin{center}
\includegraphics[width=14cm]{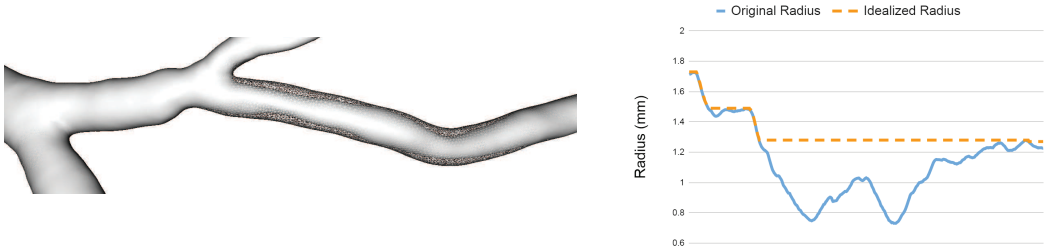}
\caption{(left) A section of a vessel with patient lumen segmentation (solid) and idealized lumen segmentation (meshed) and (right) the corresponding patient radius (blue) and idealized patient radius (dotted orange).\label{patient_ideal}}
\end{center}
\end{figure}

One of the novelties of our work is to build a patient-specific exploration space. To define the bounds on the geometric domain, at each location in the coronary tree, we consider modifications wherein the vessel radius falls between the patient radius and the idealized radius. The idealized model represents a model of a patient without regions of lumen narrowing (or ``healthy looking''). The idealized model is constructed by seeking a radius profile that is monotonically non-increasing from the ostium to the leaves and is closest to the original radius profile, using a robust objective function, formulated mathematically as \\

\begin{equation}
r_{\text{ideal}} = \text{arg min}_{r^*}  \hspace{0.2in} \Sigma_i{\sqrt{r_i^* - r_{\text{orig}, i}}} 
\end{equation}
such that
\begin{equation}
r_i^* \le r_j \text{                   if i is distal to j}
\end{equation}
\begin{equation}
r_i^* \ge m_i (\text{optional})
\end{equation}
where $r_{\text{orig},i}$ is the input radius, $r_{\text{ideal},i}$ is the idealized radius estimate and $m_i$ is an optional minimal radius constraint at a vessel tree location i. For instance, $m_i$ can be set equal to $r_{\text{orig},i}$ to force the ideal model to be point-wise at least as large as the patient model. Idealized radius values  $r_{\text{ideal},i}$ are optimized globally, over the entire coronary tree, through a dynamic programming approach. Given idealized radius values, a patient-specific idealized 3D model is constructed by dilating the original model from its local radius  $r_{\text{orig},i}$ to the target idealized radius $r_{\text{ideal},i}$. This step is performed by merging 3D spheres of radii  $r_{\text{ideal},i}$ into the implicit representation of the original geometry (signed distance field) and generating a new mesh from that dilated implicit representation. Figure~\ref{patient_ideal} shows an ideal model (and radius profile) superimposed over the patient model. We omit the index $i$ to represent the entire lumen segmentation of the appropriate radii.\par

\subsection{Flow dependent hemodynamic resistance}

Lumen geometry influences the flowrate through the model. Therefore, quantifying the impact of geometry on hemodynamics should account for the changes in the flowrate through the model. The patient-specific reduced order model (PSROM) assumes an affine resistance model for flow through the coronary arteries,
\begin{equation}
R_{\text{geom}} = a(r) + b(r)Q  
\end{equation} 
where $R_{\text{geom}}$ is the local resistance  of the lumen segmentation to blood flow, $Q$ is the blood flow rate, $a(r)$ is the resistance to flow that is purely geometric, and $b(r)$ is the sensitivity of resistance to flowrate. This  is a generalization of having a viscous (Poiseuille) and inertial (Bernoulli) term, but the coefficients being derived from HFHC simulations instead of assuming a circular tube cross-section and a parabolic velocity profile~\cite{TaylorSankaranCMAME}~\cite{UQBio}. \par

To accurately build a response surface for $a(r)$ and $b(r)$, it is imperative that we explore the space of flowrates in the patient-specific model. Since flowrate is not an independent variable, but rather depends on the boundary conditions and the vessel lumen geometry, we modify the boundary conditions to change the state of the flowrate to derive $b(r)$. We parameterize the lumen geometry by its radius, which is commonly used for analytical 1-dimensional modeling of the Navier-Stokes equations~\cite{TaylorSankaranCMAME}. This modification can be done for both the patient-specific geometry and the idealized lumen geometry, and interpolated for intermediate configurations. \par

Since we have chosen to use an affine model, two flowrates should suffice to explore the space of flowrates. Since the patient-specific model is already solved for the hyperemic state, we need only one additional flowrate. In this paper, we are interested only in dilation of the lumen geometry (not constriction) which would result only in the increase in the rate of blood flow. Therefore, we consider an additional state, called henceforth as ``superemia", wherein we reduce the microvascular resistance by 40\%. This number was chosen based on a grid search on a development dataset. Smaller values such as 10\% or 20\% tend to have less signal to noise (low flow difference and pressure loss compared to the numerical noise of the CFD simulations). Larger values such as 80\% tend to have flowrates very different from the modified lumen geometries. An illustration of the hyperemia and superemia configurations on patient and idealized geometrieis can be seen in Figure~\ref{illustration}.\par

\subsection{Patient-specific Reduced order model}

The patient-specific reduced order model (PSROM) is constructed by first calculating $a(r)$ and $b(r)$ for the original and idealized models, where two simulations are performed for each model (hyperemia and superemia). For any new configuration, the intercept can be calculated as
\begin{equation}
a(r) = \alpha a(r_{\text{orig}}) + (1 - \alpha) a(r_{\text{ideal}})
\end{equation}
where
\begin{equation}
\alpha = \frac{\alpha_{0}}{r^4} + \alpha_1
\end{equation}
\begin{equation}
\alpha_0 = \frac{r_{\text{ideal}}^4 r_{\text{orig}}^4 }{r_{\text{ideal}}^4 - r_{\text{orig}}^4}
\end{equation}
and 
\begin{equation}
\alpha_1 = \frac{r_{\text{orig}}^4}{r_{\text{ideal}}^4 - r_{\text{orig}}^4}.
\end{equation}

\begin{figure}[h!]
\begin{center}
\includegraphics[width=15cm]{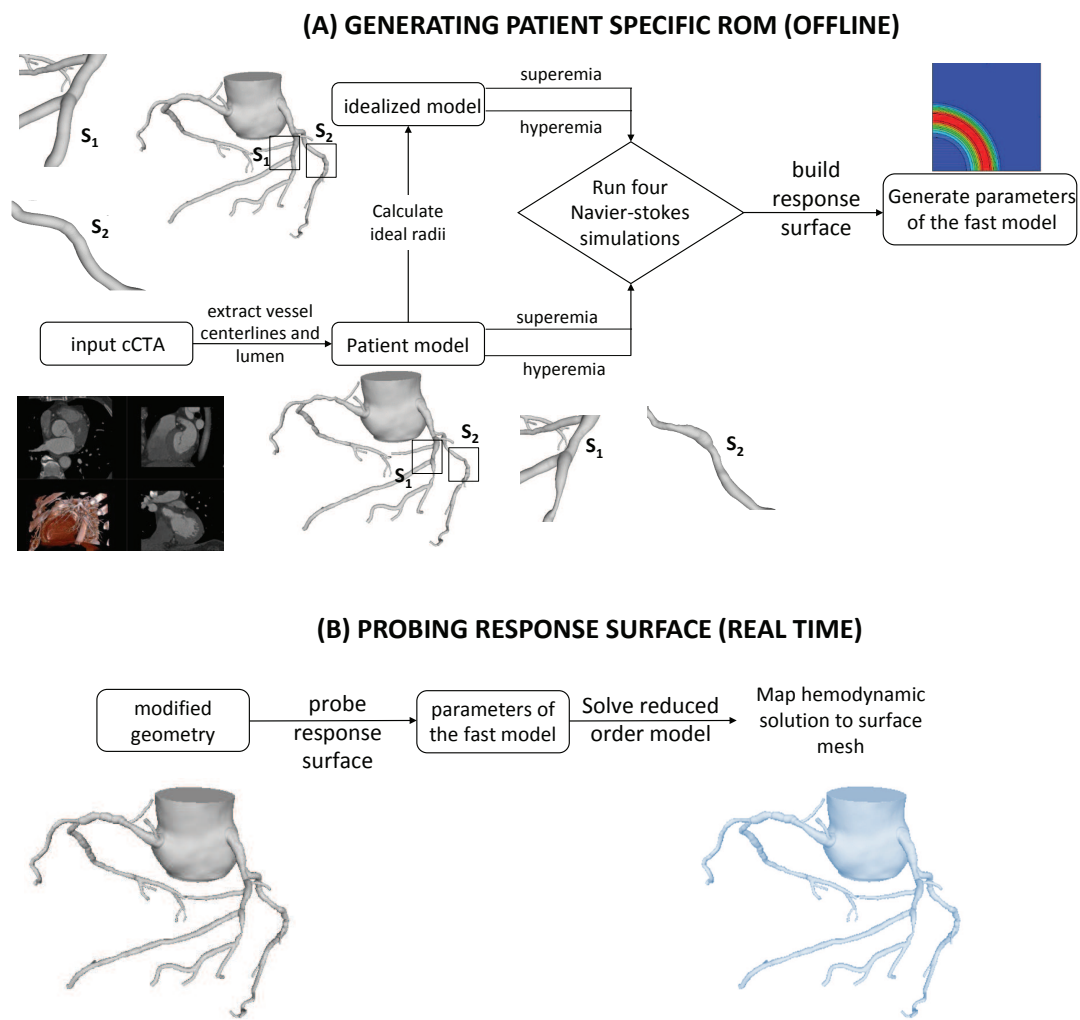}
\caption{(A) Steps involved in generating the patient specific reduced order model. First, a lumen segmentation of the patient anatomy is constructed from a coronary CT scan. Next, an idealized model of the patient's  coronary tree is constructed that represents a segmentation without lumen narrowing. Two sections of the model are highlighted to illustrate differences between the patient and idealized  models. Next, four CFD simulations are performed, two each on the patient and idealized models with two different boundary conditions. (B) the patient specific reduced order model (PSROM) can be used on any modified geometry that lies between the patient and idealized geometry to calculate the hemodynamics.\label{overall}
}
\end{center}
\end{figure}

In addition, the slope can be calculated as 
\begin{equation}
b(r) = \beta b(r_{\text{orig}}) + (1 - \beta) b(r_{\text{idealized}}),
\end{equation}
\begin{equation}
\beta = \frac{\beta_0}{A^3} \frac{\partial A}{\partial z} + \beta_1,
\end{equation}
\begin{equation}
\beta_0 = \frac{1}{\frac{1}{A_{\text{orig}}^3}\frac{\partial A_{\text{orig}}}{\partial z} - \frac{1}{(A_{\text{ideal}})^3} \frac{\partial A_{\text{ideal}}}{\partial z}}
\end{equation}
and 
\begin{equation}
\beta_1 = \frac{1}{\left(1 - \left(\frac{A_{\text{ideal}}}{A_{\text{orig}}}\right)^3 \frac{\frac{\partial A_{\text{orig}}}{\partial z}}{\frac{\partial A_{\text{ideal}}}{\partial z}}\right) }.
\end{equation}
$A_{\text{orig}}$ and $A_{\text{ideal}}$ are the patient and idealized lumen areas, and derived directly from the radii. The interpolation functions, $\alpha$ and $\beta$ were chosen to be consistent with the Poiseuille pressure loss and Bernoulli pressure loss respectively. A predictor corrector algorithm~\cite{PredictorCorrector} is used to solve the equations governing PSROM. Further details of how these quantities are calculated are given in Algorithm~\ref{algo}. Figure~\ref{overall} shows a schematic of the proposed patient-specific ROM approach. \par

\subsection{Predictor corrector algorithm}

The goal of the proposed method is to predict hemodynamics given a geometry modified from the original patient's geometry. For the modified geometry, the response surface that is constructed from the HFHC simulations is used to estimate  the coefficients $a(r)$ and $b(r)$ given a radius $r$. These are then used to calculate the local resistances $R_\text{geom}$ given a flowrate $Q$. These resistances are integrated from the leaves of the coronary tree using a circuit analogy (resistances are solved in series within a given segment, and in parallel at each bifurcation). The resistances at the ostium are then used to calculate an updated flowrate. These steps are repeated until a pre-set convergence criterion is met. The boundary conditions depend on the hemodynamics, and we account for this dependence by modeling the average size of the downstream vasculature being proportional to the pressure at the outlet (passive response to increased pressure). At each outlet, the boundary conditions are scaled from the original value by a ratio of the outlet pressures.\par

There are two conditions where we do not apply the interpolation approach. In regions just downstream of the region of lumen modification, the original patient-specific geometry may have negative resistance (increase in pressure) and thereby a negative slope. Such regions are absent when the initial section with lumen narrowing is replaced by an idealized section. We empirically define this potential pressure recovery region to be up to $20\text{mm}$ distal to a region of vessel narrowing. This was also chosen based on the development dataset. We might miss capturing pressure increase if we choose too small a pressure recovery region length cutoff. On the other hand,  the choice of larger region lengths may result in inclusion of recovery that are not related to the flow features associated with a stenosis. Further, we assume that pressure can recover only if the area gradient is negative. We therefore replace any region with negative area gradient in the pressure recovery region with a Poiseulle resistance model. In addition, we do not use the slope and intercept estimated from the four solutions if the hyperemic and superemic flowrates are too close to each other. This is because numerical convergence related errors in the HFHC model can have a non-linear effect on the slope and the intercept, and the signal to noise effect can result in lower accuracy of the prediction, especially when flowrates of the modified geometry are significantly larger than the original patient-specific model. \par 

\subsection{Overall  Algorithm}

The overall algorithm is given in Table~\ref{algo} and illustrated in Figure~\ref{overall}. A convergence criteria of $\text{tol}_2 = 2\%$ for the ostial flowrate is used. The extraction of radii from the patient-specific lumen segmentation and the subsequent generation of the idealized model are performed prior to building the physics driven response surface. ${\mathcal C}$ represents the set of all centerline points, ${\mathcal M}$ represents the set of forward neighbors (implying leaf centerline points do not have a neighbor), $P_i$ is the pressure at centerline point $i$, $Q_{i,k}$ is the flowrate of the  section connecting i to k, $a_{i,k}$ and $b_{i,k}$ are the slope and intercept of the section connecting centerline $i$ to $k$, superscripts $s$ and $h$ represent superemia and hyperemia respectively. For the predictor corrector method, the superscripts denote the iteration number. $R_{\text{eff}}$ denotes the net effective downstream resistance, and the predictor corrector method predicts the net effective resistance at each centerline point (using reverse breadth first search from leaves to ostium), corrects the ostial flow (depth first search  from  ostium to  leaves), and continues until they are consistent. For the centerline point corresponding to the outlets, the resistances are set equal to the boundary conditions, and for the rest of the centerline points, the resistances correspond to the geometric resistance~\cite{TaylorSankaranCMAME}~\cite{UQBio}. The flow-split correction factor, $\gamma_i$ is the deviation of the assumption (flow is inversely proportional to the net effective downstream resisttance) in the original patient-specific geometry, defined at each centerline point and extracted from the patient-specific hyperemic simulations. The threshold  for switching from a CFD-derived approach to a reduced order model, $\text{tol}_1$ is chosen to be $0.1$.\par

\begin{algorithm}
  \caption{Overview of the algorithm to build and predict hemodynamics using a patient-specific reduced order model.}
  \vspace{0.1in}
   \begin{center} {\sc BUILDING PHYSICS DRIVEN RESPONSE SURFACE} \end{center}
  \begin{algorithmic} 
    \State Solve the HFHC model in four configurations.
        \For{$c_i \in {\mathcal C}$}
             \For{$c_k : c_k \in {\mathcal M}(c_i)$}
                     \State $b_{i,k} = \frac{(P_i^s - P_k^s) / Q_{i,k}^s - (P_i^h - P_k^h) / Q_{i,k}^h }{Q_{i,k}^s - Q_{i,k}^h}$ \\
                     \State $a_{i,k} = \frac{(P_i^s - P_k^s)Q_{i,k}^h / Q_{i,k}^s - (P_i^h - P_k^h)Q_{i,k}^s / Q_{i,k}^h}{Q_{i,k}^h - Q_{i,k}^s}$
                     \If{$Q_{i,k}^s -  Q_{i,k}^h <= \text{tol}_1 * Q_{i,k}^h$}
                         \State $a_{i,k}=\frac{8 \mu L}{\pi r^4}$ \\
                         \State $b_{i,k}=\rho \frac{1}{A^3} \left(\frac{\partial A}{\partial z}\right)$
                     \EndIf 
            \EndFor
        \EndFor  
  \end{algorithmic}
  \vspace{-0.1in}
   \begin{center} {\sc SOLVING PSROM} \end{center}
 \vspace{-0.1in}
 
\begin{algorithmic}
     \State  Initialize n =  0, ostial flow ($\text{O}_q^0$), the superscript denoting the iteration number, $n$
        \While{ $n == 0$ or $\text{O}_q^{n} - \text{O}_q^{n-1} > \text{tol}_2 \times \text{O}_q^{n-1}$}
            \State n = n +  1
            \State $R_{i,k}^n =  a_{i,k} +  b_{i,k} Q_{i,k}^{n-1}$
            \State enqueue(leaf)
            \While ( queue not empty)
                   \State rbfs(pop queue)
           \EndWhile
            \State $\text{O}_q^n = \frac{P_{\text{aorta}}}{\text{O}_{R_{\text{eff}}}^n}$; 
             $Q^{n} = \text{dfs}(\text{ostium})$
        \EndWhile

 \Function {dfs}{current\_node = m}
    \If {is\_outlet(current\_node)}
        return
    \EndIf
     
    \If {is\_branch(current\_node)}
        $Q_{m_{d1}} = Q_m \times \frac{R_{\text{eff},m}}{R_{\text{eff},d1}} \times \gamma_{m_{d1}}$; 
        $Q_{m_{d2}} = Q_m \times \frac{R_{\text{eff},m}}{R_{\text{eff},d2}} \times \gamma_{m_{d2}}$ 
        \State dfs($m_{d1}$)
        \State dfs($m_{d2}$)
    \Else
        \For{$k : k \in {\mathcal M}(m)$}
            $Q_k = Q_m$.
            \State dfs(k)
        \EndFor
     \EndIf
     
  \EndFunction

 \Function {rbfs}{current\_node = m}
    \If {is\_ostium(current\_node)}
        return
    \EndIf
     
    \If {is\_branch(current\_node)}
        \For{$k : k \in {\mathcal M}(m)$}
            $R_{\text{local}} += 1/R_{i,k}$.
        \State $R_{\text{eff},i} = \frac{1}{R_{\text{local}}}$   
        \EndFor
    \Else
        \For{$k : k \in {\mathcal M}(m)$}
            \State $R_{\text{eff},i} += R_{i,k}$.
        \EndFor
    \EndIf
    \State enqueue(m-1)
\EndFunction
\end{algorithmic}\label{algo}
\end{algorithm}

\section{Results}
In this section, we demonstrate the accuracy and speed of the proposed method against the HFHC computational fluid dynamics model. The detailed protocol for generating the ground truth is provided in the following section. \par

\subsection{Validation protocol}
First, we chose a cohort of patients. We identified lesions in the coronary tree of each patient, and characterized them as single (focal), ostial, bifurcation or serial lesion. Focal lesions are single isolated lesions present in a vessel path from ostium to the leaves. Bifurcation lesions go across a branch, while ostial lesions occur at the junction of aorta and coronary arteries. Three or more lesions in a vessel path were characterized as serial lesions. Any lesion type that is not characterized as one of the four listed above were not considered as candidates for the validation study. \par

For each patient and lesion type, we performed a virtual remodeling of the lesion and replaced the geometry within the lesion with the idealized geometry. The rest of the patient's coronary artery tree was then blended with the remodeled region to yield a new coronary segmentation for the patient. We then performed a CFD simulation on this geometry which serves as the ground truth. The PSROM was run on the modified geometry and then compared against the ground truth simulation values. \par

In order to obtain sufficient statistical power, we chose thirteen hundred and forty one (1341) patients. Ten of these failed to yield a ground truth $\text{FFR}_{\text{CT}}$ due to meshing or convergence issues. Therefore, a total of thirteen hundred and thirty one (1331) cases were available and were used to evaluate the PSROM. \par

\begin{figure}[h!]
\begin{center}
\includegraphics[width=16cm]{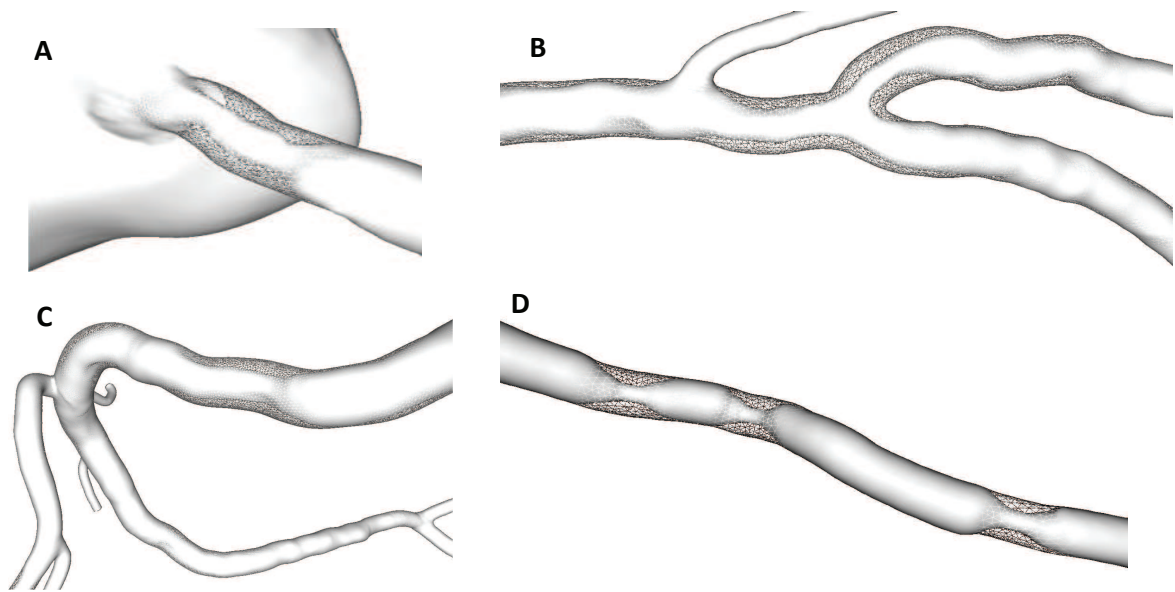}
\caption{Patient lumen segmentation (solid) in a section of the vessel superimposed with the ideal lumen segmentation (meshed) for (A) ostial lesion, (B) bifurcation lesion, (C) focal (single) lesion, and (D) serial lesion.\label{difflumen}}
\end{center}
\end{figure}

The cases were gathered from centers in the US, Japan, EU, and UK and were not used for development or validation of the algorithm. Lesions that were modified had to exhibit a lumen narrowing of at least 30\% using the metric defined in Sankaran et al.~\cite{SankaranIEEE} with $r_{\text{ideal}}$ being used for the healthy or ideal radius. Samples of the different lesion types along with their idealized lumen geometry are shown in Figure~\ref{difflumen}. \par

For each case, we evaluated results at salient points downstream of the lumen modification but avoided bifurcation or pressure recovery regions. In some cases, more than one comparison location was identified. In total, 2578 comparisons across 1331 patient-specific models were  performed. \par

In order to demonstrate equivalence, we tested that the bias is within $0.01$ and the standard deviation is within $0.02$.  The numbers are motivated primarily by a mesh independence study on the ground truth values, which demonstrated that if we were to replace $\text{FFR}_{\text{CT}}$ with a very fine mesh, the 95\% confidence interval of error is within 0.03, which translates to a standard deviation of $\sim 0.02$. To demonstrate equivalence in bias, we used a two one-sided t-test (TOST) and to demonstrate equivalence in standard deviation, we used a chi-squared test. The null hypothesis for the TOST test is the absolute difference between the proposed method and ground truth is less than $0.03$. The null hypothesis for the chi-squared test is that the standard deviation is greater than 0.02. \par

\subsection{Performance}

\begin{figure}[h]
\begin{center}
\includegraphics[width=15cm]{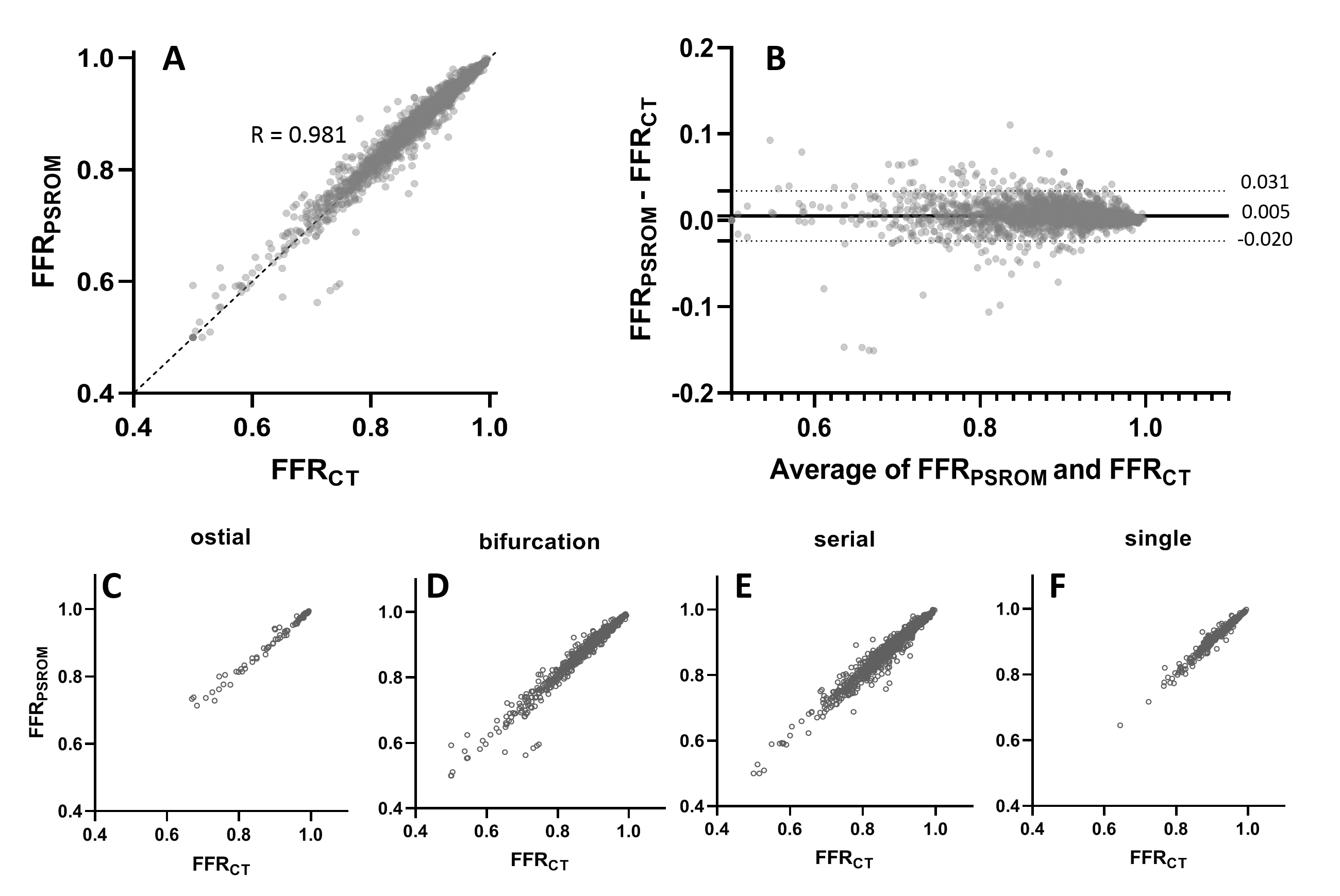}
\caption{(A) Scatterplot comparing the ground truth ($\text{FFR}_{\text{CT}}$) solution obtained by solving the Navier-Stokes equations (HFHC model) to the solution obtained using the patient-specific reduced order model ($\text{FFR}_{\text{PSROM}}$). (B) Bland-Altman plot showing the error as a function of the average values, along with the bias  and 95\% limits of agreement, (C-F) Scatterplot comparing $\text{FFR}_{\text{PSROM}}$ to $\text{FFR}_{\text{CT}}$ in modifications with ostial, bifurcation, serial and single (focal) disease respectively.\label{scatter}}
\end{center}
\end{figure}

The overall bias in the evaluated dataset was 0.00524, with a corresponding p-value using TOST (two one-sided t-test) of $<$1e-3. The standard deviation in the evaluated dataset was 0.0147, with a corresponding p-value using chi-squared test of $<$1e-3. The correlation coefficient between $\text{FFR}_{\text{CT}}$ and PSROM was 0.982 (95\% CI 0.981-0.984, p$<$0.001). The 95\% Bland-Altman limits of agreement between $\text{FFR}_{\text{CT}}$ and $\text{FFR}_{\text{PSROM}}$ was (-0.020, 0.031) as seen in Figure~\ref{scatter}. These are within the 95\% confidence intervals of measurement reproducibility. Slope and intercept were observed to be 0.972 and 0.030 respectively. \par

The mean difference of $\text{FFR}_{\text{CT}}$ before and after modification was 0.111 with a standard deviation of 0.121. The maximum change in $\text{FFR}_{\text{CT}}$ before and after lumen modification was 0.70. For those locations where the change in $\text{FFR}_{\text{PSROM}}$ before and after modification was greater than 0.1 (N = 1082), the bias and standard deviation were  0.006 and 0.020 respectively. For locations where the change in FFR using PSROM before and after modification was less than 0.1 (N = 1496), a bias of 0.004 and standard deviation of error of 0.009 was observed. This also demonstrates that the proposed  method is able to predict $\text{FFR}_{\text{CT}}$ for large and small changes before and after modification. Ninety five percent of all the comparisons had an error of less than 0.03 compared to $\text{FFR}_{\text{CT}}$.\par 

\begin{figure}[h!]
\begin{center}
\includegraphics[width=17cm]{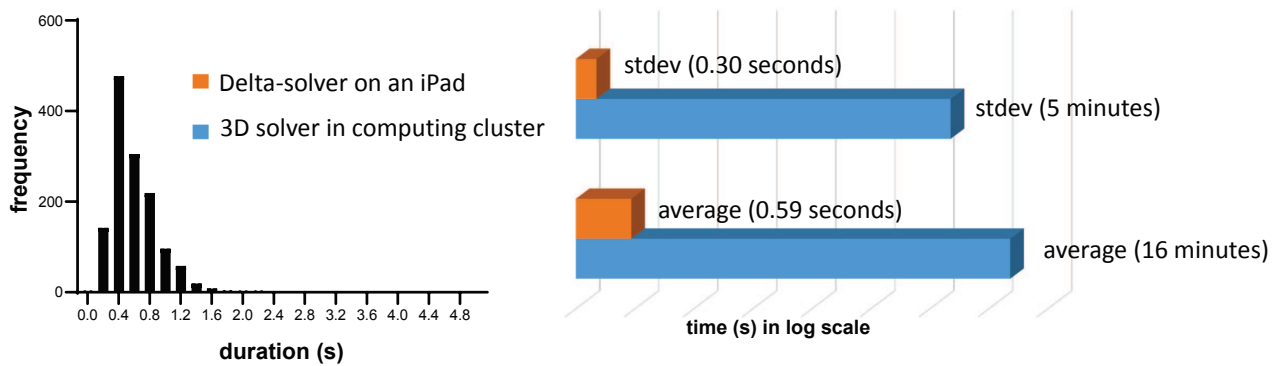}
\caption{(left) Run times of the proposed method with a median of 0.53 seconds, (right) A comparison of the average and standard deviation (stdev)
 of the full solver and the PSROM method plotted in a logarithmic scale. \label{DStiming} }
\end{center}
\end{figure}

The average time taken to execute the PSROM algorithm for the prediction step was $0.59 \pm 0.30$ sec (median 0.53 sec). The maximum time taken was 2.22 sec with a 75th percentile value of 0.75 sec. Figure~\ref{DStiming} shows the histogram of duration of the PSROM algorithm on an iPad Pro 10.5 device. Average and standard deviation of run times across the entire population are shown for the ground truth and the PSROM methods (in logarithmic scale).\par

\begin{table}
\centering
\begin{tabular}{ |c | c | c | c | c |}
\hline
lesion type & sample size & bias & standard deviation & 95\% CI \\
\hline
\hline
bifurcation & 782 & 0.0051 & 0.0171 & (-0.028, 0.039) \\
ostial & 78 & 0.0089 & 0.0153 & (-0.021, 0.039) \\
focal & 465 & 0.0030 & 0.0149 & (-0.014, 0.019) \\
serial & 1253 & 0.0059 & 0.0084 & (-0.023, 0.035) \\
overall & 2578 & 0.0052 & 0.0148 & (-0.020, 0.031) \\
\hline
\end{tabular}
\caption{Performance of the proposed method categorized by lesion type. Bias, standard deviation and 95\% confidence intervals are reported along with the sample size.}
\label{stratifyperf}
\end{table}

\begin{table}
\centering
\begin{tabular}{ |c | c | c | c |}
\hline
range & N & bias & standard deviation \\
\hline
\hline
$[0.00,0.70)$ & 74 & 0.0160 & 0.0261 \\
$[0.70,0.75)$ & 84  & 0.0035 & 0.0394 \\
$[0.75,0.80)$ & 151 & 0.0081 & 0.0214 \\
$[0.80,0.85)$ & 309 & 0.0069 & 0.0165 \\
$[0.85,0.90)$ & 558 & 0.0068 & 0.0149 \\
$[0.90,1.00]$ & 1402 & 0.0035 & 0.0080 \\
\hline
\end{tabular}
\caption{Performance of the proposed method stratified by $\text{FFR}_{\text{PSROM}}$. Bias and standard deviation are reported along with the sample size.}
\label{rangestratify}
\end{table}

Further, the results were analyzed based on stratification by lesion type. For bifurcation lesions (N = 782), bias was 0.005 (p$<$1e-6) and standard deviation of error was 0.017 (p$<$1e-6). For serial lesions (N = 1253), bias was 0.006 (p$<$1e-6) and standard deviation of error was 0.008 (p$<$1e-6). For focal lesions (N = 465), bias was 0.003 (p$<$1e-6) and standard deviation of error was 0.015 (p$<$1e-6). For ostial lesions (N = 78), bias was 0.009(p=0.25) and standard deviation of error was 0.015 (p=0.004). The N corresponds to the number of comparisons, which may include more than one per patient. Table~\ref{stratifyperf} summarizes these results. \par

Table~\ref{rangestratify} summarizes performance of the algorithm when data was stratified based on $\text{FFR}_\text{PSROM}$ in different ranges (or buckets). Each bucket has a different sample size which is expected based on a typical randomized patient population. The bias is the highest in the bucket with smallest $\text{FFR}_{\text{CT}}$ and the standard deviation is highest in the $[0.70,0.75)$ bucket. \par

\begin{figure}[h]
\begin{center}
\includegraphics[width=15cm]{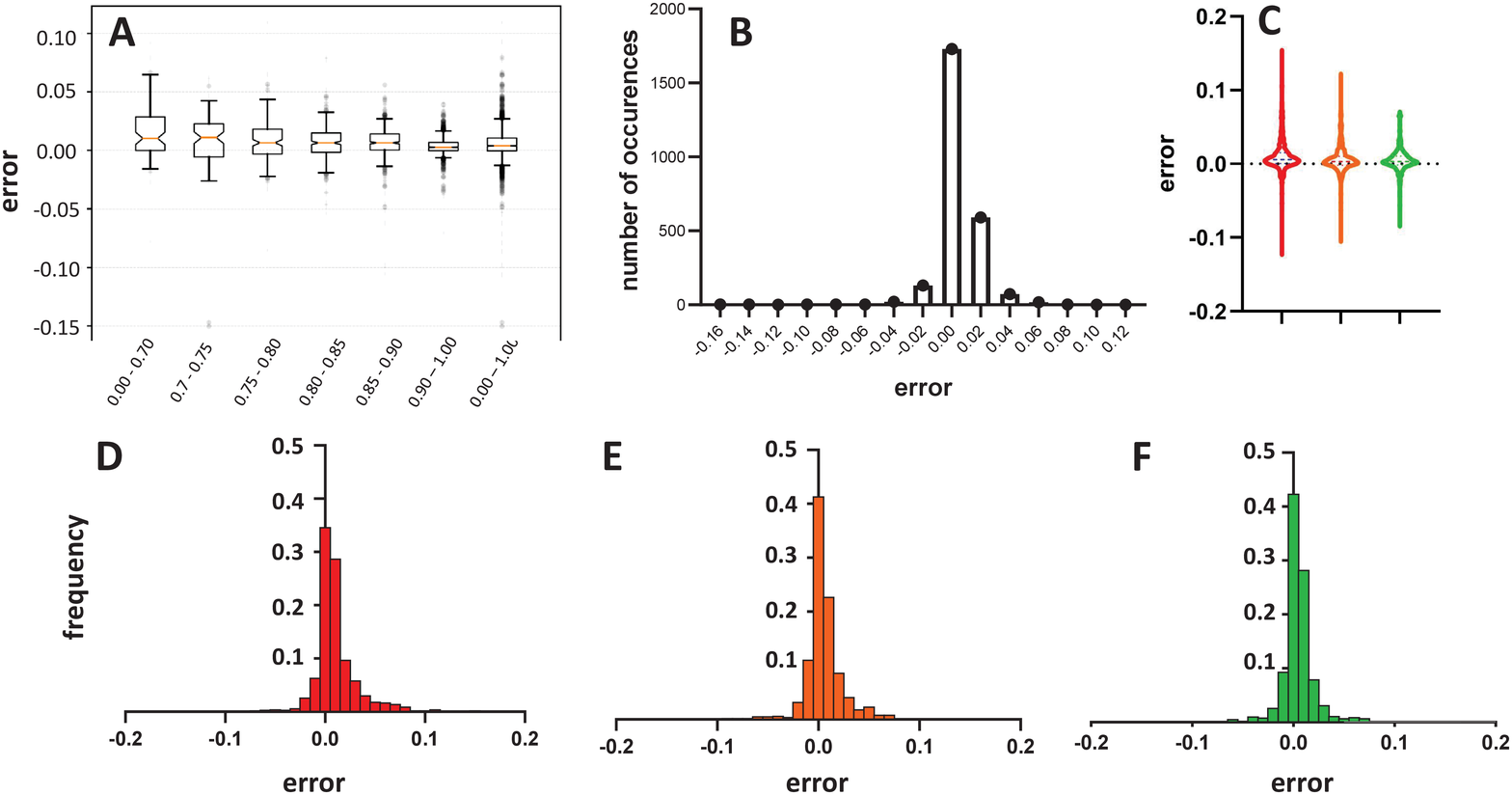}
\caption{Illustration of the properties of the error between $\text{FFR}_\text{CT}$ and $\text{FFR}_\text{PSROM}$ distributions. (A) Error distribution in different $\text{FFR}_\text{PSROM}$ ranges. (B) error histogram in the evaluation dataset, (C-F) error distributions corresponding to (red) unconstrained, (orange) hybrid based on a ROM cutoff, and (green) hybrid with a pressure recovery fix term, which was used for all the results in this paper.\label{errordist}}
\end{center}
\end{figure}

Figure~\ref{errordist} demonstrates error properties of the FFR results of the PSROM method. It demonstrates that the error distribution is centered close to a mean of zero, and has a normal distribution around the mean. In addition, the figure also demonstrates the added value of the proposed method against a ROM either not accounting for the vessels with low flow difference or not accounting for the pressure recovery. \par

\begin{figure}[h!]
\begin{center}
\includegraphics[width=15cm]{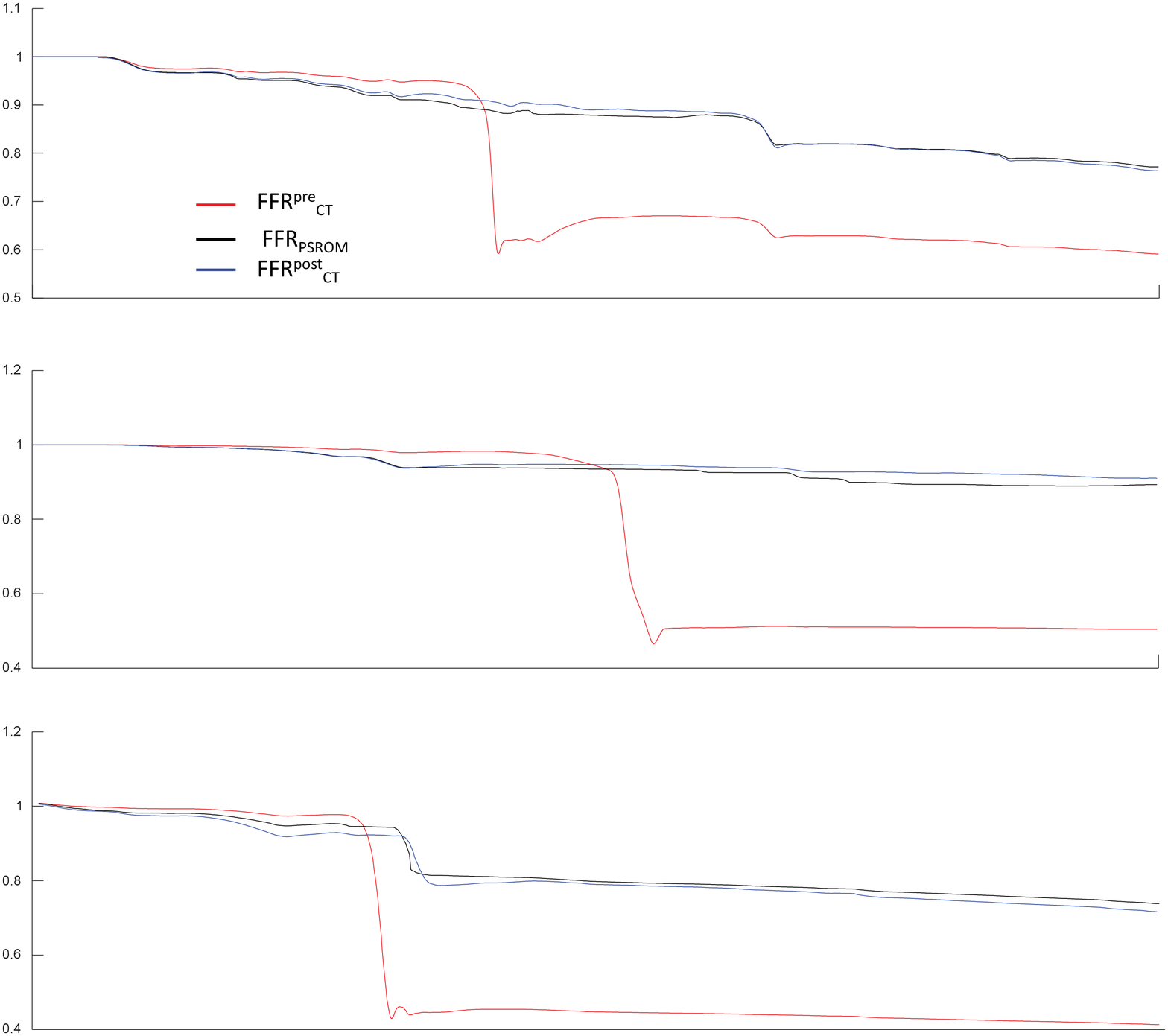}
\caption{FFR tracing from the ostium to a distal end of a modified vessel, comparing the original trace, PSROM prediction and ground truth.\label{trace}}
\end{center}
\end{figure}

Figure~\ref{trace} shows pressure tracing in three vessels with serial lesions with superimposed traces of the $\text{FFR}_{\text{CT}}$ of the patient model before  and after modification and $\text{FFR}_{\text{PSROM}}$. These examples demonstrate the ability of the PSROM method to capture the pressure drop in tandem lesions with higher flowrate post lumen modification. \par

\subsection{Validation against invasive data}

We also clinically validated the PSROM method using invasive data obtained from percutaneous coronary intervention (PCI) or stenting procedures. PCI is the process of inserting a catheter with a balloon attached to restore and attain ideal lumen geometry on a section of the blood vessel. The objective of this procedure is to relieve symptoms  of ischemia and mitigate risk of coronary artery disease. Some of the biggest challenges of PCI is to know a-priori which lesions have the largest impact on FFR (for serial lesions), how stenting may restore blood flow (for all lesion types, with diffuse atherosclerosis potentially not increasing coronary blood flow even after restoring vessel caliber), or how to size the stent (radius and length). The PSROM method was first applied to a patient with serial lesions where the invasive FFR was also measured~\cite{Ihdayhid}. Two lesions (A - proximal and B - distal) were identified in the left arterior descending (LAD) artery.  Invasive revascularization of B and a combination of A and B were performed. We applied the PSROM, blinded, to these two configurations, and demonstrated excellent match for the two scenarios (lesion B - invasive 0.77, PSROM 0.79, and lesion A + B - invasive 0.85, PSROM 0.88). We also applied the PSROM method to a cohort of thirteen patients with serial lesions who underwent PCI~\cite{Modi}. The use of PSROM improved the correlation coefficient of the comparison to 0.75 over a value of 0.44 using conventional FFRct (p $<$ 0.001). Further, the bias and coefficient of variation of PSROM was 0.01 and 7\% respectively, compared to 0.05 and 37\% using pre-PCI values. \par
\section{Discussion}

We formulated and developed a method, called PSROM, that enables fast assessment of hemodynamics in response to modifying a patient-specific geometry. The method was built on the two principles. First, resistance of lumen geometry to flow was modeled to have a linear relationship with the flowrate, with the slope and intercept being functions of arterial geometry. This was primarily motivated by the steady-state solution to a 1-dimensional Navier-Stokes equation that has a linear and quadratic term for pressure loss. Second, flow split at bifurcations was modeled to be inversely related to the net downstream resistance of the respective vessels. Net downstream resistance accounts for both the geometric resistance as well as the boundary conditions. Additionally, we assumed that the idealized model represents the largest possible dilation of the CT derived lumen geometry, though the method outlined in this work is applicable even if the strategy to define idealized model is different.\par

We initially built a local space around a patient-specific geometry that was constrained by a maximum possible geometric expansion using a patient-specific idealized model. The idealized model was derived from the patient's lumen segmentation by fitting a radius profile that is non-increasing from the root to the leaves. The idealized model represents a patient-specific lumen segmentation without regions of lumen narrowing that is the closest in L-1/2 norm with the original lumen segmentation. Further, we used a sampling of different microvascular properties (dilation) to derive the properties of the affine model. In this work, we define this space using four configurations, two different lumen segmentations (patient and idealized model) and two boundary conditions each. The two geometric configurations were chosen to be the extrema of the lumen segmentation (patient and idealized). The boundary conditions were chosen to be hyperemia, and another lower resistance that we call superemia, which has 40\% less resistance and allows more flow through the vessels. This, in turn, enables the PSROM to parameterize the intercept and slope. A physics driven response surface was used to represent the resistance of geometry to flow within this patient-specific exploration space. \par

We deviated from the formulation to consider two special situations - (i) if the vessel is severely diseased such that flowrate is limited by geometry, then there is insufficient flow separation between the hyperemic and superemic conditions. In such situations, we used a model that is not derived from the HFHC simulations, but rather uses the analytical reduced order model (1D Navier-Sokes equations) constrained by the hyperemic patient and ideal model hemodynamics, (ii) in the sub-regions within the coronary tree where pressure recovers, i.e. pressure increases as we traverse the vessel tree from the proximal to distal section. These are regions where a reduced order model interpolation does not work because geometric resistance is negative. Therefore, we replace the interpolatory behavior by using an analytical model in the ``pressure recovery" zone (that we define in this work as within 20 mm from the stent). The results showed that these developments were needed to achieve the desired performance. \par
In order to accurately capture pressure recovery, models for smaller scale fluid structures such as eddies and turbulence need to be captured. The value of 20 mm used in this paper was based on empirical observation in the development dataset, but this can be improved further by understanding the link between lumen geometry, flowrate and the resultant pressure recovery region length and extent. \par

The method we proposed was validated against a large cohort of patients with different lesion characteristics: namely, bifurcation, serial, single and ostial lesions. The method performed very well against the ground truth CFD simulations in over 2500 evaluation points. This was enabled by deriving the ROM directly from set of CFD simulations, thereby negating the scenarios where other low dimensional models perform poorly due to a combination of various factors such as (i) assuming a parabolic velocity profile, (ii) approximating geomety by radius, and/or (iii) failing to account for turbulent losses. Leveraging the patient-specific model to model pressure recovery, except the region just downstream of the modified region, also helped improve the performance over traditional reduced order models. The method also performed well against invasive data on a cohort of patients who underwent PCI. To our knowledge, no published method in literature has demonstrated such good performance. The technique was recently cleared by the U.S. Food and Drug Administration (FDA) for clinical use. \par

The proposed method showed significant discriminatory power in predicting the $\text{FFR}_{\text{CT}}$ post modification compared to the ground truth three dimensional CFD simulations. Additionally, we demonstrated a good correlation coefficient, Bland-Altman limits of agreement, and agreement across lesion types. Two sets of clinical studies~\cite{Modi,Ihdayhid}, where the invasive measurements were blinded, showed promising results when compared with PSROM. \par

The method performed well in identifying ``silent lesions'', those that are masked by another severe lesion in the same vessel path. Relieving the lumen narrowing in a flow-limiting lesion results in more flow through the tandem (less severe) lesion (refer Fig.~\ref{trace}). Consequently, this will result in a higher pressure drop across the unmodified lesion.  The results we obtained demonstrated that the proposed method is able to capture this increase in pressure drop accurately. \par

The proposed method does not fully account for changes to the patient post-intervention such as potential plaque shift, side branch jailing, and the material used for lumen dilation. Further, the changes to the microvascular resistance are modeled only with regards to passive remodeling of the vessel wall (relationship of area to flow). Finally, the performance in modifications with a large change was not on par with the performance in healthier sections. \par

We believe that the method outlined in this work has the potential to aid interventional procedural planning by its ability to virtually model different scenarios before performing the actual intervention. This could potentially increase cath lab efficiency by reducing unnecessary procedures, procedural duration, total radiation dose as well as highlight scenarios that might be missed by visual inspection of an angiogram.

 \bibliographystyle{unsrt}

\bibliography{Planner_MethodsPaper_arxiv}

\end{document}